\documentclass[a4paper,11pt,leqno]{amsart}

\usepackage{amssymb,hyperref}
\usepackage[latin1]{inputenc}
\usepackage[francais,english]{babel}
\usepackage{graphicx}
\usepackage{color}

\usepackage[T1]{fontenc}
\usepackage[latin1]{inputenc}
\usepackage[all]{xy}
\usepackage{stmaryrd}

\usepackage{psfrag}
\usepackage{epsfig}
\usepackage{wrapfig}
\usepackage{amsmath}
\usepackage{amssymb,latexsym}
\usepackage[mathscr]{eucal}
\usepackage{url}
\usepackage{fancybox}
\swapnumbers

\numberwithin{equation}{section}

\DeclareMathAlphabet{\mathrmsl}{OT1}{cmr}{m}{sl}

\newcommand{\oper}[3][n]{\newcommand{#2}{\mathop{\mathrm{#3}}%
\ifx n#1\nolimits\else\limits\fi} }
\newcommand{\rsoper}[3][n]{\newcommand{#2}{\mathop{\mathrmsl{#3}}%
\ifx n#1\nolimits\else\limits\fi} }

\newcounter{mnotecount}[section]
\renewcommand{\themnotecount}{\thesection.\arabic{mnotecount}}
\newcommand{\mnote}[1]
{\protect{\stepcounter{mnotecount}}$^{\mbox{\footnotesize  $
      \bullet$\themnotecount}}$ \marginpar{\raggedright\tiny\em
    $\!\!\!\!\!\!\,\bullet$\themnotecount: #1} }

\newcommand{\eq}[1]{\ref{#1}}

\renewcommand{\qed}{\hfill $\blacksquare$ \medskip \\}
\newcommand{\preuve}{\noindent {\sc Proof:\ }}

\theoremstyle{main}
\newtheorem{thm} {\bf  Theorem} [section]

\newtheorem{lem} [thm] {\bf  Lemma}
\newtheorem{prop} [thm] {\bf  Proposition}

\newtheorem{déf}[thm]{\bf Definition}
\newtheorem{Rq}[thm]{\bf  Remark }

\DeclareFontFamily{OT1}{rsfs}{}
\DeclareFontShape{OT1}{rsfs}{m}{n}{ <-7> rsfs5 <7-10> rsfs7 <10-> rsfs10}{}
\DeclareMathAlphabet{\mycal}{OT1}{rsfs}{m}{n}

\renewcommand{\div}{\operatorname{div}}

\newcommand{\Vol}{\operatorname{Vol}}
\newcommand{\Scal}{\operatorname{Scal}}
\newcommand{\Ric}{\operatorname{Ric}}
\newcommand{\tr}{ {\rm tr}}
\newcommand{\dtr}{ {\rm dtr}}

\newcommand{\Ker}{\operatorname{Ker}}

\newcommand{\End}{\operatorname{End}}

\newcommand{\AdS}{\operatorname{AdS}}

\renewcommand{\geq}{\geqslant}
\renewcommand{\leq}{\leqslant}
\newcommand{\eps}{\varepsilon}

\renewcommand{\Im}{\operatorname{Im}}
\newcommand{\R}{{\Bbb R}}
\newcommand{\C}{{\Bbb C}}

\newcommand{\D}{{\frak D}}
\renewcommand{\k}{{\bf K}}

\title[eigenvalue estimate for the dirac--witten operator]{Optimal eigenvalue estimate for the Dirac--Witten operator on bounded domains with smooth boundary}
\author{Daniel MAERTEN}
\date{\today}
\address{LMPT, Fédération Denis Poisson, Faculté des Sciences, Parc de Grandmont, F--37200 Tours, France}
\email{{\color{red}daniel.maerten@yahoo.fr}}

\begin{document}

\maketitle

\begin{abstract}
Eigenvalue estimate for the Dirac--Witten operator is given on bounded domains (with smooth boundary) of spacelike hypersurfaces satisfying the dominant energy condition, under four natural boundary conditions (MIT, APS , modified APS and chiral conditions). Roughly speaking, any eigenvalue of the Dirac--Witten operator satisfies
$$ \left|\lambda\right|^{2} \geq \frac{n}{(n-1)}\frak R_{0} \ ,$$
where $\frak R_{0}$ is the infinimum of (the opposite of) the lorentzian norm of the constraints vector. Equality cases are also investigated and lead to interesting geometric situations.
\end{abstract}

\section{Introduction}

The aim of spectral geometry is to derive some geometric properties from the study of the spectrum of a certain elliptic operator, which is mostly a Laplace operator (Laplacian on functions, Hodge Laplacian on $p$--forms or Dirac--Laplacian on spinors). In that context, a classical issue is to give lower bounds for the eigenvalues of the Dirac operator on a compact manifold $(M^{n},g)$, $(n\geq2)$. In \cite{F}, T.~Friedrich proved that any eigenvalue $\lambda$ of the Dirac operator satisfies
\begin{equation}\label{Friedrich}
	\lambda^{2} \geq \frac{n}{4(n-1)}\inf_{M}\Scal^{g} \ ,
\end{equation}
where $\Scal^{g}$ denotes the scalar curvature of $(M^{n},g)$. A few years later, O.~Hijazi \cite{Hij1,Hij2} improved this result by showing that, for any $n\geq3$,
\begin{equation}\label{Hijazi1}
	\lambda^{2} \geq \frac{n}{4(n-1)}\mu_{1} \ ,
\end{equation}
where $\mu_{1}$ is the first eigenvalue of the conformal Laplacian. Clearly, estimates (\eq{Friedrich}) and (\eq{Hijazi1}) are under interest only if the scalar curvature is positive. In a recent work,  O.~Hijazi and X.~Zhang \cite{HijZ} established an analogous version of (\eq{Friedrich}) for the Dirac--Witten operator under a chiral boundary condition, of a compact spacelike hypersurface $(M,g,k)$ satisfying the dominant energy condition, namely
\begin{equation}\label{Hijazi2}
\lambda^{2} \geq \frac{n}{4(n-1)}\inf_{M}\left\{\Scal^{g} + \left(\tr _{g}k \right)^{2}-\left|k\right|^{2}_{g} - 2\left|\delta_{g}k+ \dtr_{g}k \right|_{g}\right\} \ .
\end{equation}
On the other hand, S.~Raulot \cite{R} proved an eigenvalue estimate for the Dirac operator on domains with boundary. More precisely, if $\Omega$  is a compact domain of an $n$--dimensional Riemannian spin manifold $(M,g)$, whose boundary $\partial\Omega$ has positive mean curvature $H$, he showed  under a natural boundary condition (called "MIT" boundary condition), that any eigenvalue $\lambda$ of the spectrum of the Dirac operator on $\Omega$ (which is an unbounded discrete set of complex numbers of positive imaginary part) satisfies
\begin{equation}\label{Raulot}
\left|\lambda\right|^{2} \geq \frac{n}{4(n-1)}\inf_{\Omega}\Scal^{g} + n\Im(\lambda)\inf_{\partial\Omega}H \ .
\end{equation}
The most interesting fact in Inequality~(\eq{Raulot}) is that the bound depends on some boundary geometric quantity, which is not the case for Inequality~(\eq{Hijazi2}). In addition, S.~Raulot showed that equality in (\eq{Raulot}) leads to the existence  of imaginary Killing spinor on  $\Omega$, and also to the conclusion that the boundary $\partial\Omega$ is a totally umbilical and constant mean curvature hypersurface.\\

The goal of this article is to generalise Inequalities~(\eq{Hijazi2}) and (\eq{Raulot}) in several directions. Indeed, we shall prove analogous versions of (\eq{Hijazi2}) for the Dirac--Witten operator on bounded domains, under four natural boundary conditions (see \cite{HijMR}) and also generalise (\eq{Raulot}) for the Dirac--Witten operator.\\
The article is organised as follows: In Section~\eq{notations}, we give our geometric conventions and preliminary results; Section~\eq{preuves} is devoted to the statements of the main results and their proves.


\section{Geometric Background}\label{notations}

\subsection{Notations}

We consider $(N^{n+1},\gamma )$ a Lorentzian manifold of signature $(-,+, \cdots,+)$ which contains  $(M^{n},g,k)$,  a spin (in dimension 3 this only means orientable) Riemannian hypersurface (not necessarily compact) whose induced metric is $g$  and second fundamental form (extrinsic curvature) is $k$. Let $\Omega$  be a compact domain in  $(M^{n},g,k)$ satisfying the dominant energy condition, which reads as the following inequality along $\Omega$
$$\Scal^{g} + \left(\tr_{g} k\right)^{2}-\left|k\right|^{2}_{g} \geq2\left|\delta_{g}k+ \dtr_{g}k \right|_{g} \ .$$
We will work with the complex spinor bundle of $N$ restricted to the hypersurface domain $\Omega$, that is to say $\Sigma :=\Sigma (N)_{|\Omega}$ which is given by the choice of a unit normal $e_{0}$ of $M$ in $N$, along $\Omega$. More precisely, if one denotes by $P_{\text{Spin}(n,1)}(N)$ the bundle of $\text{Spin}(n,1)$--frames on $N$, and by $\rho_{n,1}$ the standard representation of  $\text{Spin}(n,1)$ then
$$\Sigma(N):= P_{\text{Spin}(n,1)}(N)\times_{\rho_{n,1}}\Bbb C^{[(n+1)/2]}.$$
Now the choice of unit normal $e_{0}$ of $M$ in $N$, along $\Omega$, induces a natural inclusion ${\text{Spin}(n)\subset\text{Spin}(n,1)}$ and so we can define
$$\Sigma := P_{\text{Spin}(n,1)}(N)_{|\Omega}\times_{(\rho_{n,1})_{|\text{Spin}(n)}}\Bbb C^{[(n+1)/2]}.$$
$\Sigma$ naturally carries two sesquilinear inner products: the first one denoted by $(\ast ,\ast  )$  is $\text{Spin}(n,1)$-invariant (it is not necessary positive); the second one which is denoted by ${\left\langle\ast ,\ast \right\rangle:=(e_{0}\cdot \ast ,\ast )}$  is $\text{Spin}(n)$-invariant and Hermitian definite positive ($\cdot $ is the Clifford action with respect to the metric $\gamma $). The Hermitian or anti-Hermitian character of the Clifford multiplication by vectors differs if we consider $(\ast ,\ast  )$ or $\left\langle\ast ,\ast \right\rangle$ and is described in \cite{M},  for instance. $\Sigma $ is also endowed with two different connections $\nabla, \overline{\nabla}$ which are respectively the Levi-Civita connections of $\gamma$ and $g$. Let us take a spinor field $\psi \in \Gamma(\Sigma)$ and a vector field $X\in\Gamma(T\Omega)$, then our conventions are
\begin{eqnarray*}
\nabla_{X}\psi & = &\overline{\nabla}_{X}\psi - \frac{1}{2}k(X)\cdot e_{0} \cdot \psi \\
\left\langle k(X),Y \right\rangle_{\gamma} &= & \left\langle \nabla_{X}Y,e_{0}\right\rangle_{\gamma} \quad .
\end{eqnarray*}
In these formulae $\cdot$ denotes the Clifford action with respect to the metric $\gamma$. The induced metric, Levi-Civita connection and second fundamental form of the boundary $\partial \Omega$, are respectively denoted by $\ell,\widetilde{\nabla},\theta $. Our conventions are, for any $X,Y \in \Gamma(T\partial\Omega)$ and $\psi\in \Sigma_{|\partial\Omega}$,
\begin{eqnarray*}
  \overline{\nabla}_{X}Y&=&\widetilde{\nabla}_{X}Y +\theta (X,Y)\nu\\
  \overline{\nabla}_{X}\psi&=&\widetilde{\nabla}_{X}\psi+\frac{1}{2}\theta (X)\cdot\nu\cdot\psi,\\
\end{eqnarray*}
where $\nu$ is the unit normal to $\partial \Omega$ pointing  inside $\Omega$, and $\cdot$ still denotes the Clifford action with respect to the metric $\gamma$. Finally we define the following  geometric quantity
\begin{equation*}
	\frak R _{0}:= \frac{1}{4}\inf_{\Omega}\left\{\Scal^{g} + \left(\tr _{g}k \right)^{2}-\left|k\right|^{2}_{g} - 2\left|\delta_{g}k+ \dtr_{g}k \right|_{g}\right\} \ .
\end{equation*}


\subsection{Dirac-Witten operators and Bochner--Lichnerowicz formul\ae}

From now on, $\left(e_{k}\right)^{n}_{k=0}$ denotes an orthonormal basis at the point, with respect to the metric $\gamma$, and where $\nu= e_{1}$. We can then define the Dirac--Witten operator of $\nabla$
$$ \D\varphi= \sum^{n}_{k=1}e_{k}\cdot \nabla_{e_{k}}\varphi \  ,$$
where $\cdot$ is the Clifford action with respect to the metric $\gamma$. Notice that $\D$ can be considered as a deformation of $\overline{\D}$ the usual Dirac operator of $\overline{\nabla}$  since we recover $\D=\overline{\D}$ as soon as $k\equiv 0$. The Dirac-Witten operator $\mathfrak{D}$ is clearly formally self adjoint in $L^{2}$ with respect to $\left\langle \ast,\ast\right\rangle$ in the class of compactly supported spinor fields, and we have the classical Bochner-Lichnerowicz-Weitzenböck formula (cf.~\cite{Bartnik, H1, PT} for instance) 
\begin{equation}\label{Bochner1}
	\D^{*}\D= \D\D= \nabla^{*} \nabla + \frak R  \ ,
\end{equation}
where 
\begin{eqnarray*}
	\frak R&:=& \frac{1}{4} \big(\text{Scal}^{\gamma} + 4 \text{Ric}^{\gamma}(e_{0},e_{0}) + 2e_{0}\cdot \text{Ric}^{\gamma}(e_{0}) \big)\\
	&=& \frac{1}{4}\left\{   \big(\Scal^{g} + \left(\tr _{g}k \right)^{2}-\left|k\right|^{2}_{g}\big) + 2(\delta_{g}k+ \dtr_{g}k ) \cdot e_{0}  \right\} \ .
\end{eqnarray*}
As usual, we derive from (\eq{Bochner1}) an integration formula. The idea is to consider a certain spinor field $\varphi$ and to define the 1--form ${\omega_{\varphi}\in\Gamma(T^{*}\Omega)}$ by the relation
$$ \omega_{\varphi}(X)= \left\langle \nabla_{X}\varphi+X\cdot \D\varphi,\varphi\right\rangle \ .$$
Then,  computing the $g$--divergence of $\omega_{\varphi}$, we obtain
$$\div_{g}\omega_{\varphi}= \left|\D \varphi\right|^{2}-\left|\nabla \varphi\right|^{2}-\left\langle \frak R \varphi,\varphi \right\rangle\ ,$$
which gives, applying Stokes' theorem
\begin{equation}\label{Int1}
	\int_{\Omega} \left|\D \varphi\right|^{2}= \int_{\Omega} \left|\nabla \varphi\right|^{2} + \int_{\Omega} \left\langle \frak R \varphi,\varphi \right\rangle + \int_{\partial\Omega}\left\langle \nabla_{\nu}\varphi +  \nu\cdot\D \varphi,\varphi \right\rangle \ .
\end{equation}

We now have to introduce $P$ the twistor operator with respect to the connection $\nabla$, which is defined by the relation
$$P_{X} \varphi:= \nabla_{X}\varphi + \frac{1}{n}X\cdot \D \varphi \ ,$$
for every $X\in\Gamma(T\Omega)$ and every spinor field $\varphi\in\Gamma(\Sigma)$. We derive a second integration formula based on a Bochner--Lichnerowicz formula for the twistor operator.

\begin{prop}
For any spinor field $\varphi\in\Gamma(\Sigma)$, we have
\begin{equation}\label{Int2}
	\int_{\Omega} \left|P \varphi\right|^{2}= \left(\frac{n-1}{n}\right)\int_{\Omega} \left|\D \varphi\right|^{2}  - \int_{\Omega} \left\langle \frak R \varphi,\varphi \right\rangle -\int_{\partial\Omega}\omega_{\varphi}(\nu) \ .
\end{equation}
\end{prop}

\preuve We first prove the Bochner--Lichnerowicz formula for $P$, namely
\begin{eqnarray*}
	\left|P\varphi\right|^{2}&=& \sum^{n}_{k=1} \left\langle \nabla_{e_{k}}\varphi + \frac{1}{n}e_{k}\cdot \D \varphi, \nabla_{e_{k}}\varphi + \frac{1}{n}e_{k}\cdot \D \varphi\right\rangle \\
	&=& \left|\nabla\varphi\right|^{2} + \left(\frac{1}{n}-\frac{2}{n}\right)\left|\D\varphi\right|^{2} \\
	&=& \left|\nabla\varphi\right|^{2}-\frac{1}{n}\left|\D\varphi\right|^{2} \ .
\end{eqnarray*}
We integrate this formula on $\Omega$ so that
\begin{eqnarray*}
	\int_{\Omega} \left|P \varphi\right|^{2}&=&  \int_{\Omega} \left|\nabla\varphi\right|^{2} - \frac{1}{n}\int_{\Omega} \left|\D \varphi\right|^{2} \\
	&=& \left(\frac{n-1}{n}\right)\int_{\Omega} \left|\D \varphi\right|^{2}  - \int_{\Omega} \left\langle \frak R \varphi,\varphi \right\rangle -\int_{\partial\Omega}\omega_{\varphi}(\nu) \ ,
\end{eqnarray*}
where the second line is obtained thanks to (\eq{Int1}). \qed
It is clear that the value of the 1--form $\omega_{\varphi}$ depends upon the boundary condition that will be used. We will need the intermediate result:
\begin{lem} 
For any spinor field $\varphi\in\Gamma(\Sigma)$, we have along the boundary $\partial\Omega$
\begin{equation}\label{integrand}
	\omega_{\varphi}(\nu)=\left\langle \nu\cdot  e_{j}\cdot \widetilde{\nabla}_{e_{j}}\varphi + \frac{1}{2}\Big( (\tr_{g} k ) \nu\cdot e_{0}\cdot - k(\nu)\cdot e_{0}\cdot + (\tr_{\ell}\theta)  \Big) \varphi,\varphi\right\rangle \ .
\end{equation}
\end{lem}

\preuve Just compute, using the relations of compatibility between the different connections
\begin{eqnarray*}
	\nabla_{\nu}\varphi+\nu\cdot \D \varphi &=& \nu\cdot \sum^{n}_{j=2} e_{j}\cdot \nabla_{e_{j}}\varphi \\
	&=& \nu\cdot  e_{j}\cdot \left(\widetilde{\nabla}_{e_{j}}-\frac{1}{2}k(e_{j})\cdot e_{0}\cdot + \frac{1}{2}\theta(e_{j})\cdot\nu\cdot\right)\varphi  \\
	&=& \nu\cdot  e_{j}\cdot \widetilde{\nabla}_{e_{j}}\varphi + \frac{1}{2}\Big( (\tr_{g} k ) \nu\cdot e_{0}\cdot - k(\nu)\cdot e_{0}\cdot + (\tr_{\ell}\theta)  \Big)\varphi \ .
\end{eqnarray*}
\qed


\section{Main Results}\label{preuves}

\subsection{"MIT" boundary condition}
Remind that our aim is to find an estimate for the spectrum of the Dirac--Witten operator under a natural boundary condition (that had been used in order to obtain some \textit{black hole} version of the positive mass theorem for asymptotically hyperbolic manifolds \cite{ChH,M}). It consists on finding a lower bound for  $\left|\lambda\right|^{2}$ where $\lambda$ is any non--zero complex (a priori) number involving in the following elliptic first order boundary problem
$$(\textrm{MIT}) \quad \left\{
\begin{array}{ll}
\D \varphi= \lambda\varphi & \textrm{ on } \Omega \\
F(\varphi)=  \varphi & \textrm{ on } \partial\Omega
\end{array} \right. \ ,
$$
where $\varphi$ is non--zero (eigen--)spinor field, and where $F\in\End\left(\Sigma_{|\partial\Omega}\right)$ is defined  by the relation  ${F(\psi)=i\nu\cdot\psi}$. The boundary condition $F(\psi)=  \psi$ was originally introduced by physicists of the MIT, and this the reason why it is  often called "MIT" boundary condition. Their idea is based on the fact that the Dirac (and also  the Dirac--Witten) operator on manifold with boundary is not formally self--adjoint anymore, since we have the following integration by parts formula
\begin{equation}\label{Int3}
	\int_{\Omega}\left\langle \D \varphi,\psi\right\rangle =\int_{\Omega}\left\langle \varphi,\D \psi\right\rangle - \int_{\partial\Omega}\left\langle \nu\cdot\varphi,\psi\right\rangle \ ,
\end{equation}
for any spinor fields $\varphi,\psi$. This defect of self--adjointness has a consequence on the spectrum of the problem (MIT).
\begin{lem}
The spectrum of {\rm (MIT)} is a discrete set of complex numbers with positive imaginary parts. 
\end{lem}

\preuve Let $\lambda$ be any eigenvalue of ${\rm (MIT)}$ with $\varphi$ a corresponding eigenspinor field. Just take $\psi=i\varphi$ in (\eq{Int3}) and consequently get
$$ 2\Im(\lambda)\int_{\Omega} \left|\varphi\right|^{2}= \int_{\partial\Omega} \left|\varphi\right|^{2} \geq 0 \ ,$$
and so $\Im(\lambda)\geq0$. Assume now that $\Im(\lambda)=0$, then $\varphi$ should vanish on $\partial\Omega$ and so on the whole $\Omega$, by the continuation principle. This contradicts the fact that an eigenspinor is by definition non identically zero. \qed

For later use we introduce some geometric quantities.
\begin{déf} We set
$$H_{0}^{{\rm MIT}}:= \inf_{\partial\Omega}\left\{\tr_{\ell}\theta - \left|k^{\partial\Omega}(\nu)\right|_{g}\right\} \  , \ k^{\partial\Omega}(\nu):=\sum^{n}_{j=2}k(\nu,e_{j})e_{j} \ .$$
\end{déf}

The first main result of this note is a generalisation of the lower bound of \cite{R} for the Dirac--Witten operator (we recover (\eq{Raulot}) when we set $k\equiv 0$ on $\Omega$).

\begin{thm}\label{mainMIT}
Let $\Omega$ be a compact domain of a spacelike spin hypersurface $(M,g,k)$ which satisfies the dominant energy condition along $\Omega$ (so that $\frak R _{0}\geq0$). The boundary $\partial\Omega$ is assumed to verify $H^{{\rm MIT}}_{0}>0$. Then under the {\rm (MIT)} boundary condition, the spectrum of the Dirac--Witten operator on $\Omega$ is an unbounded discrete set of complex numbers with positive imaginary part, such that  any eigenvalue satisfies
\begin{equation}\label{spectreMIT}
	\left|\lambda\right|^{2} \geq \frac{n}{(n-1)}\Big(\frak R_{0} +  H_{0}^{{\rm MIT}}\Im(\lambda)\Big) \ .
\end{equation}
\end{thm}

\preuve
The fact that, under the (MIT) boundary condition, the spectrum of the Dirac--Witten operator on $\Omega$ is an unbounded discrete set of complex numbers with positive imaginary part, has been proved in the previous lemma.\\
As far as Inequality~(\eq{spectreMIT}) is concerned, we  consider a 1--parameter family of modified spinorial Levi--Civita connection. Indeed, for any $\alpha\in\R$, we define the action of  $\nabla^{\alpha}$ on $\Sigma$ by the relation
$$ \nabla^{\alpha}_{X}\varphi:= \nabla_{X}\varphi +i \alpha X\cdot \varphi \ ,$$
for every $X\in\Gamma(T\Omega)$.  We can then define the Dirac--Witten operator with respect to the Killing connection $\nabla^{\alpha}$
$$\D^{\alpha}\varphi= \sum^{n}_{k=1}e_{k}\cdot \nabla^{\alpha}_{e_{k}}\varphi \ ,$$
where $\cdot$ is the Clifford action with respect to the metric $\gamma$. An easy computation gives the relation $\D^{\alpha}= \D -i n \alpha$, so that $\D^{\alpha}$ is not formally self--adjoint in $L^{2}$ with respect to $\left\langle \ast,\ast\right\rangle$ in the class of compactly supported spinor fields if $\alpha\neq 0$.\\
We consider a certain spinor field $\varphi$ and  define the 1--form ${\omega^{\alpha}_{\varphi}\in\Gamma(T^{*}\Omega)}$ by the relation
$$ \omega^{\alpha}_{\varphi}(X)= \left\langle \nabla^{\alpha}_{X}\varphi+X\cdot \D^{\alpha}\varphi,\varphi\right\rangle \ .$$
We notice that $\omega^{\alpha}_{\varphi}(X)= \omega_{\varphi}(X) - \alpha(n-1)\left\langle i X\cdot\varphi,\varphi \right\rangle$.
Then,  we only have to compute  the $g$--divergence of the 1--form $\xi(X):=\left\langle i X\cdot\varphi,\varphi \right\rangle$. To this end, we assume that our local base satisfies $\overline{\nabla}_{e_{j}}e_{m}=0$ at the point where the computation is made (this is equivalent to $\nabla_{e_{j}}e_{m}=-k(e_{j},e_{m})e_{0}$) so that 
\begin{eqnarray*}
	\div_{g}\xi &=& -\sum^{n}_{j=1} e_{j}\cdot \xi(e_{j}) \\
	&=& - \left\langle i \underbrace{\overline{\nabla}_{e_{j}}e_{j}}_{=0}\cdot \varphi,\varphi \right\rangle - \left\langle i e_{j}\cdot \overline{\nabla}_{e_{j}}\varphi, \varphi\right\rangle - \left\langle i e_{j}\cdot \varphi, \overline{\nabla} _{e_{j}}\varphi \right\rangle \\
	&=& 	 - \left\langle i \D\varphi,\varphi \right\rangle - \frac{1}{2}\left\langle i \underbrace{e_{j}\cdot k(e_{j})}_{=-\tr_{g}k} \cdot e_{0} \cdot \varphi, \varphi\right\rangle  - \left\langle \varphi, i\D \varphi \right\rangle + \frac{1}{2}\left\langle i \underbrace{k(e_{j}) \cdot e_{j}}_{=-\tr_{g}k} \cdot e_{0} \cdot \varphi, \varphi\right\rangle \\
	&=&  - \left\langle i \D\varphi,\varphi \right\rangle - \left\langle \varphi, i\D \varphi \right\rangle  \ .
\end{eqnarray*}
Using Stokes' theorem and integration formula~(\eq{Int1}), we get
$$ \int_{\partial\Omega }  \omega^{\alpha}_{\varphi}(\nu) = \int_{\Omega} \left|\D \varphi\right|^{2} - \int_{\Omega} \left|\nabla \varphi\right|^{2} -\int_{\Omega} \left\langle \frak R \varphi,\varphi,\right\rangle  + \alpha(n-1) \int_{\Omega} \big( \left\langle i\D \varphi,\varphi\right\rangle + \left\langle \varphi,i \D \varphi \right\rangle\big) \ .$$ 
 Easy computations lead to the relations
\begin{eqnarray*}
	\left|\nabla^{\alpha}\varphi\right|^{2} &=& \left|\nabla\varphi\right|^{2} +n\alpha^{2}\left|\varphi\right|^{2} +\alpha\big(\left\langle i\D\varphi,\varphi\right\rangle + \left\langle \varphi, i\D \varphi\right\rangle\big) \\
	\left|\D^{\alpha}\varphi\right|^{2}&=& \left|\D \varphi\right|^{2} +  n^{2}\alpha^{2}\left|\varphi\right|^{2} + n\alpha \big(\left\langle i\D\varphi,\varphi\right\rangle + \left\langle \varphi, i\D \varphi\right\rangle\big) \ .
	\end{eqnarray*}
Plugging this into our integration by parts formula, we obtain
\begin{equation}\label{Int4}
\int_{\partial\Omega }  \omega^{\alpha}_{\varphi}(\nu) = \int_{\Omega} \left|\D^{\alpha} \varphi\right|^{2} - \int_{\Omega} \left|\nabla^{\alpha} \varphi\right|^{2} -\int_{\Omega} \left\langle \frak R^{\alpha} \varphi,\varphi,\right\rangle   \ ,
\end{equation}
where $\frak R^{\alpha}:= \frak R +n(n-1)\alpha^{2}$. The twistor operator with respect to the connection $\nabla^{\alpha}$ is defined by 
$$P^{\alpha}_{X} \varphi:= \nabla^{\alpha}_{X}\varphi + \frac{1}{n}X\cdot \D^{\alpha} \varphi \ ,$$
for every $X\in\Gamma(T\Omega)$ and every spinor field $\varphi\in\Gamma(\Sigma)$. We straightly have
\begin{equation*}
	\left|P^{\alpha} \varphi \right|^{2} = \left|\nabla^{\alpha}\varphi\right|^{2} - \frac{1}{n}\left|\D^{\alpha}\varphi\right|^{2} \geq 0 \ ,
\end{equation*}
and, after integration on $\Omega$ we get, using (\eq{Int4})
\begin{equation}\label{Int5}
	0 \leq \int_{\Omega}\left|P^{\alpha} \varphi \right|^{2} = \left(\frac{n-1}{n}\right) \int_{\Omega} \left|\D^{\alpha} \varphi\right|^{2}  -\int_{\Omega} \left\langle \frak R^{\alpha} \varphi,\varphi,\right\rangle  - \int_{\partial\Omega }  \omega^{\alpha}_{\varphi}(\nu) \ .
\end{equation}

We  now need to give some elementary properties of the boundary spinorial endomorphism $F$ (the proof is left to the reader).

\begin{prop}\label{symMIT}
The endomorphism $F$ is symmetric, isometric with respect to $\left\langle\ast,\ast\right\rangle$, commutes to the action of $\nu\cdot$  and anticommutes to each $e_{k}\cdot,\, (k \neq 1)$.
\end{prop}
The important fact is that the boundary condition of (MIT) allows us to express the boundary integrand of (\eq{Int5}) in terms of the extrinsic curvature tensors $k$ and $\theta$. 

\begin{lem}\label{F}
If $F(\varphi)=\varphi$ then along the boundary $\partial\Omega$ we have
\begin{equation}\label{integrandMIT}
	\omega^{\alpha}_{\varphi}(\nu)=\frac{1}{2}\left\langle e_{0}\cdot \left\{ \big(\tr_{\ell} \theta - 2(n-1)\alpha\big)e_{0} + k^{\partial\Omega}(\nu) \right\} \cdot\varphi,\varphi\right\rangle  \ .
\end{equation}
\end{lem}
\noindent
\textsc{Proof of Lemma~\eq{F}:}
We use Equation~(\eq{integrand}) 
\begin{eqnarray*}
	\omega^{\alpha}_{\varphi}(\nu) &=& \left\langle \nu\cdot  e_{j}\cdot \widetilde{\nabla}_{e_{j}}\varphi+ \alpha i \nu\cdot e_{j}\cdot e_{j}\cdot\varphi + \frac{1}{2}\Big( (\tr_{g} k ) \nu\cdot e_{0}\cdot - k(\nu)\cdot e_{0}\cdot + (\tr_{\ell}\theta)  \Big) \varphi,\varphi\right\rangle\\
	&=& \left\langle \nu\cdot  e_{j}\cdot \widetilde{\nabla}_{e_{j}}\varphi- (n-1) \alpha i \nu\cdot \varphi + \frac{1}{2}\Big( (\tr_{\ell} k ) \nu\cdot e_{0}\cdot - k^{\partial\Omega}(\nu)\cdot e_{0}\cdot + (\tr_{\ell}\theta)  \Big) \varphi,\varphi\right\rangle \ .
\end{eqnarray*}
Taking Proposition~\eq{symMIT} into account and also the assumption $F(\varphi)=\varphi$, it comes out
\begin{eqnarray*}
	\omega^{\alpha}_{\varphi}(\nu) &=& \frac{1}{2}\left\langle \big(  (\tr_{\ell}\theta -2\alpha(n-1))- k^{\partial\Omega}(\nu)\cdot e_{0}\cdot\big)\varphi,\varphi  \right\rangle \\
	&=& \frac{1}{2}\left\langle e_{0}\cdot \big(  (\tr_{\ell}\theta -2\alpha(n-1))e_{0}+ k^{\partial\Omega}(\nu)\big)\cdot\varphi,\varphi  \right\rangle \\
	&=& \frac{1}{2}\left\langle e_{0}\cdot \k^{\alpha}_{{\rm MIT}}	 \cdot \varphi, \varphi  \right\rangle \ ,
\end{eqnarray*}
where we have defined the vector field $ \k^{\alpha}_{{\rm MIT}} \in \Gamma(TN_{|\partial\Omega})$ as
$$ \k^{\alpha}_{{\rm MIT}} = \big(\tr_{\ell}\theta -2\alpha(n-1)\big)e_{0}+ k^{\partial\Omega}(\nu) \ .$$\qed
Then, it is well known that the boundary integrand $\omega^{\alpha}_{\varphi}(\nu)$ is non--negative if and only if $ \k^{\alpha}_{{\rm MIT}} $ is causal and future oriented, which reads as
\begin{equation}\label{bordMIT}
	\tr_{\ell}\theta -2\alpha(n-1) \geq \left|k^{\partial\Omega}(\nu)\right|_{g} \ .
\end{equation}
As a consequence, if one assumes that (\eq{bordMIT}) holds and that $\D\varphi=\lambda\varphi$, then (\eq{Int5}) implies
\begin{eqnarray*}
	0 &\leq&  \left(\frac{n-1}{n}\right) \left|\lambda-in\alpha \right|^{2}\int_{\Omega} \left| \varphi\right|^{2}  -\int_{\Omega} \left\langle \frak R^{\alpha} \varphi,\varphi,\right\rangle \\
	0 &\leq&  \left(\frac{n-1}{n}\right)\left(\left|\lambda\right|^{2} + n^{2}\alpha^{2} -2n\alpha\Im(\lambda)\right)\int_{\Omega} \left| \varphi\right|^{2} -\int_{\Omega} \left\langle \frak R^{\alpha} \varphi,\varphi,\right\rangle \\
	0 &\leq&   \left(\frac{n-1}{n}\right)\left(\left|\lambda\right|^{2} -2n\alpha\Im(\lambda)\right)\int_{\Omega} \left| \varphi\right|^{2} -\int_{\Omega} \left\langle \frak R \varphi,\varphi,\right\rangle \ ,
\end{eqnarray*}
which is possible only if $ \left|\lambda\right|^{2} \geq   \frac{n}{(n-1)}\frak R_{0} +2n\alpha\Im(\lambda)\ $. Now just take $\alpha=\alpha_{0}= \frac{1}{2(n-1)}H^{{\rm MIT}}_{0}$ so that $\k^{\alpha_{0}}_{{\rm MIT}}$ is clearly causal and future oriented and we obtain the desired inequality, namely
$$\left|\lambda\right|^{2} \geq   \frac{n}{(n-1)}\Big(\frak R_{0} + H^{{\rm MIT}}_{0}\Im(\lambda)\Big)\ .$$~\qed

The equality case in (\eq{spectreMIT}) is not easy to treat in general, but we can however deduce some information under a natural additional assumption.

\begin{thm}\label{MITegal}
Under the assumptions of Theorem~\eq{mainMIT}, equality in (\eq{spectreMIT}) always leads to the existence of an imaginary $\nabla$--Killing spinor on $\Omega$ of Killing number $-i \alpha_{0}=-\lambda/n$. If we assume furthermore that ${k^{\partial\Omega}(\nu)=0}$ along the boundary $\partial\Omega$, then  $\partial\Omega$ is a totally umbilical and constant mean curvature hypersurface of $\Omega$. 
\end{thm}

\preuve
Suppose now that equality holds in~(\eq{spectreMIT}). Thereby, there exists a non--zero $\lambda$--eigenspinor field $\varphi$ such that 
\begin{equation}\label{egalMIT}
	P^{\alpha_{0}}\varphi \equiv 0 \ , \quad \left\langle e_{0}\cdot\k_{{\rm MIT}}^{\alpha_{0}}\cdot\varphi,\varphi \right\rangle \equiv 0 \ .
\end{equation}
The second equation of~(\eq{egalMIT}) simply says that the vector field $\k_{{\rm MIT}}^{\alpha_{0}}$ is lightlike, which means in other words
$$\tr_{\ell}\theta - \left|k^{\partial\Omega}(\nu)\right|\equiv 2(n-1)\alpha_{0} \ .$$
The first equation of~(\eq{egalMIT}) can be reformulated as 
$$ \nabla_{X}\varphi + \left(\frac{\lambda}{n}\right) X\cdot \varphi =0 \ , \quad\ \forall X \in\Gamma(T\Omega),$$
that is, $\varphi$ is a Killing spinor. It is well known that the Killing number $\left(\frac{-\lambda}{n}\right)$ has to be either real or purely imaginary. But remind that $\Im(\lambda)>0$, and so $\lambda \in i \R^{*}_{+}$. For later use, we set $\lambda=i\mu$ with $\mu$ a positive real number.
We consider now any vector field $X\in\Gamma(T\partial\Omega)$ tangent to the boundary, and we compute
\begin{eqnarray*}
	 -i\left(\frac{\mu}{n}\right) X \cdot \varphi &=& \nabla_{X}\varphi\\
	&=& \nabla_{X}(i\nu\cdot\varphi) \\
	&=& i \big(\nabla_{X}\nu\cdot\varphi + \nu\cdot \nabla_{X}\varphi\big)\\
	&=& i \left(  \overline{\nabla}_{X}\nu\cdot -k(X,\nu)e_{0}\cdot -i \left(\frac{\mu}{n}\right)\nu\cdot X\cdot  \right)\varphi \\
	&=& i \left( -\theta(X)\cdot -k(X,\nu)e_{0}\cdot  + \left(\frac{\mu}{n}\right) X\cdot  \right)\varphi \ .
\end{eqnarray*}
This finally leads to
\begin{equation}\label{meanMIT1}
	\forall X\in\Gamma(T\partial\Omega) \ , \quad \left( -\theta(X)\cdot -k(X,\nu)e_{0}\cdot  + \left(\frac{2\mu}{n}\right) X\cdot  \right)\varphi=0 \ .
\end{equation}
The condition~(\eq{meanMIT1})  notably implies that for any  vector field $X\in\Gamma(T\partial\Omega)$, the vector field $\left\{-\theta(X) -k(X,\nu)e_{0}  + \left(\frac{2\mu}{n}\right)X\right\} $ is lightlike, and consequently we have 
\begin{equation}\label{meanMIT2}
	\forall X\in\Gamma(T\partial\Omega) \ , \quad k(X,\nu)^{2}=\left| \theta(X) -2\left(\frac{\mu}{n}\right) X  \right|^{2}_{g} \ .
\end{equation}
The striking relation~(\eq{meanMIT2}) somehow measures the failure of $\partial\Omega$ to be totally umbilical. We now prove that $\mu=n\alpha_{0}$. Indeed, using (\eq{meanMIT1}), we compute along the boundary 
\begin{eqnarray*}
	\left(\frac{2\mu}{n}\right) \sum^{n}_{j=2} \left\langle e_{j}\cdot \varphi, e_{j}\cdot\varphi\right\rangle &=&  \frac{2(n-1)\mu}{n}\left|\varphi\right|^{2}\\
	&=& \sum^{n}_{j=2} \left\langle e_{j}\cdot \varphi,\left\{ \theta(e_{j})\cdot +k(e_{j},\nu)e_{0}\right\}\cdot\varphi\right\rangle\\
	&=&-\sum^{n}_{j=2}\left\langle  \varphi,e_{j}\cdot\left\{ \theta(e_{j})\cdot +k(e_{j},\nu)e_{0}\right\}\cdot\varphi\right\rangle\\
	&=& \left\langle \varphi, \left\{\tr_{\ell} \theta - k^{\partial\Omega}(\nu)\cdot e_{0}\cdot\right\} \varphi\right\rangle\\
	&=&  \left\langle \varphi, e_{0}\cdot\k^{0}_{{\rm MIT}}\cdot\varphi \right\rangle \ ,
\end{eqnarray*}
which can be written in short as $\left\langle \varphi, e_{0} \cdot \k^{0}_{{\rm MIT}}\varphi \right\rangle=\frac{2\mu(n-1)}{n}\left|\varphi\right|^{2}$. Besides,  using the second equation of (\eq{egalMIT}), we know that 
\begin{eqnarray*}
	0&=& \int_{\partial\Omega}\left\langle  e_{0}\cdot\k^{\alpha_{0}}_{{\rm MIT}}\cdot\varphi,\varphi \right\rangle\\
	&=& \int_{\partial\Omega}\left\langle \varphi, e_{0}\cdot\k^{0}_{{\rm MIT}}\cdot\varphi \right\rangle - 2(n-1)\alpha_{0}\int_{\partial\Omega}\left|\varphi\right|^{2}\\
	&=& 2(n-1)\left(  \frac{\mu}{n} - \alpha_{0}\right) \int_{\partial\Omega}\left|\varphi\right|^{2} \ ,
\end{eqnarray*}
which entails $\mu=n\alpha_{0}$ since the eigen--spinor $\varphi$ cannot vanish identically along $\partial\Omega$. Suppose furthermore, that $k^{\partial\Omega}(\nu)=0$ on the boundary, then it is clear from (\eq{meanMIT2}) that $\partial\Omega$ is a totally umbilical and constant mean curvature hypersurface. \qed


\subsection{APS Boundary Condition}
The Atiyah--Patodi--Singer (APS) boundary condition was often used in order to prove positive mass theorem for asymptotically flat black holes (that is to say for asymptotically flat manifolds with boundary \cite{H2, M1}). More precisely, we denote by $\Pi_{\pm }$ the $L^{2}$-orthogonal projections on the spaces of eigenspinors of positive (respectively negative) eigenvalues of the Dirac operator of the boundary $\widetilde{\D}:=\underset{j\geq2}{\sum}e_{j}\cdot_{\partial\Omega} \widetilde{\nabla}_{e_{j}}$, where $\cdot_{\partial\Omega}$ denotes the Clifford action with the respect to boundary metric $\ell$. In this section, our aim is to find  a lower bound for  $\left|\lambda\right|^{2}$ where $\lambda$ is any non--zero complex (a priori) number involving in the following elliptic first order boundary problem
$$(\textrm{APS}) \quad \left\{
\begin{array}{ll}
\D \varphi= \lambda\varphi & \textrm{ on } \Omega \\
\Pi_{+}(\varphi)=  0 & \textrm{ on } \partial\Omega
\end{array} \right. \ ,
$$
where $\D$ still denotes the Dirac--Witten operator. We first prove that the spectrum of (APS) is real.
\begin{lem}
The spectrum of {\rm (APS)} is a discrete set of real numbers. 
\end{lem}

\preuve Let $\lambda\in \C$ be any eigenvalue of ${\rm (APS)}$ with $\varphi$ a corresponding eigenspinor field. Now,  $\cdot_{N},\ \cdot_{\Omega} ,\ \cdot_{\partial\Omega}$ the Clifford actions of respectively $\gamma, \ g, \  \ell$ satisfy the following relations
\begin{equation}\label{clifford}
	 X \cdot _{\partial \Omega} \psi = X\cdot _{\Omega}\nu \cdot _{\Omega}\psi,  \ X \cdot _{ \Omega} \psi = i X\cdot _{N}e_{0} \cdot _{N}\psi \ .
\end{equation}
These identities imply  $\widetilde{\D}=-\nu\cdot \underset{j\geq2}{\sum}e_{j}\cdot \widetilde{\nabla}_{e_{j}}$, where $\cdot = \cdot _{N}$ is still the Clifford action of $\gamma$. The spinor bundle $\Sigma_{|\partial\Omega}$ has a natural $L^{2}$--orthogonal decomposition

\begin{equation}\label{decomposition}
\begin{array}{ccccccc}
	\Sigma_{|\partial\Omega}&=&\Im \Pi_{+} & \oplus & \Im \Pi_{-} &\oplus &\Ker \widetilde{\D} \\
	&=& \Sigma^{+}_{|\partial\Omega}&\oplus&\Sigma^{-}_{|\partial\Omega}&\oplus & \Sigma_{|\partial\Omega}^{0}
\end{array} \ .
\end{equation}
An important fact is that $\widetilde{\D}$ anticommutes with $\nu\cdot$, namely $\nu\cdot\widetilde{\D}=- \widetilde{\D}\nu\cdot $. Indeed, for any spinor $\psi\in \Sigma_{|\partial\Omega}$ we have
\begin{eqnarray*}
	\widetilde{\D}(\nu\cdot \psi )&=& -\nu\cdot  \sum^{n}_{j=2} e_{j}\cdot \widetilde{\nabla}_{e_{j}}(\nu\cdot\psi)\\
	&=& -\nu\cdot e_{j}\cdot  \left( \overline{\nabla}_{e_{j}}(\nu\cdot\psi) - \frac{1}{2}\theta(e_{j})\cdot\nu\cdot\nu\cdot\psi     \right) \\
	&=&  -\nu\cdot  e_{j}\cdot  \left(-\theta(e_{j})\cdot\psi + \nu\cdot \overline{\nabla}_{e_{j}}\psi + \frac{1}{2}\theta(e_{j})\cdot\psi \right)   \\
	&=&  -\nu\cdot  e_{j}\cdot  \left( - \frac{1}{2} \nu\cdot\theta(e_{j})\cdot\nu\cdot\psi + \nu\cdot \overline{\nabla}_{e_{j}}\psi  \right) \\
	&=& - \nu\cdot \widetilde{\D}\psi \ .
\end{eqnarray*}
Thereby, it turns out that $\nu\cdot \Sigma^{\pm}_{|\partial\Omega}\subset\Sigma^{\mp}_{|\partial\Omega}$ and $\nu\cdot\Sigma^{0}_{|\partial\Omega}\subset \Sigma^{0}_{|\partial\Omega}$. Remind our integration by parts formula (\eq{Int3}) 
\begin{equation*}
	\int_{\Omega}\left\langle \D \varphi,\psi\right\rangle =\int_{\Omega}\left\langle \varphi,\D \psi\right\rangle - \int_{\partial\Omega}\left\langle \nu\cdot\varphi,\psi\right\rangle \ .
\end{equation*}
Taking $\varphi$ a solution of (APS) for an non--zero eigenvalue $\lambda$ and setting $\psi=\varphi$ in the formula above, leads to $$(\lambda-\overline{\lambda})\left\|\varphi\right\|_{L^{2}(\Omega)}=-\left\langle \nu\cdot\varphi,\varphi\right\rangle_{L^{2}(\partial\Omega)} \ .$$ But, $\Pi_{+}\varphi=0$ means $\varphi_{|\partial\Omega} \in \Sigma^{-}_{|\partial\Omega}\oplus \Sigma^{0}_{|\partial\Omega}$, and $\nu\cdot \Sigma^{\pm}_{|\partial\Omega}\subset\Sigma^{\mp}_{|\partial\Omega} + $ gives $\left\langle \nu\cdot\varphi,\varphi\right\rangle_{L^{2}(\partial\Omega)}=0$ because of decomposition~(\eq{decomposition}). We thus obtain $\lambda-\overline{\lambda}=0$, that is $\lambda$ is real. \qed

For later use we introduce a geometric quantity.
\begin{déf} We set
$$H_{0}^{{\rm APS}}:= \inf_{\partial\Omega}\left\{\tr_{\ell}\theta - \left((\tr_{\ell}k)^{2} + \sum^{n}_{j=2} k(\nu,e_{j})^{2} \right)^{\frac{1}{2}}  \right\}  \ .$$
\end{déf}

The second main result of this note is a generalisation of the lower bound of \cite{HijMR} for the Dirac--Witten operator under the (APS) boundary condition (we recover the inequality of \cite{HijMR} when we set $k\equiv 0$ on $\Omega$).

\begin{thm}\label{mainAPS}
Let $\Omega$ be a compact domain of a spacelike spin hypersurface $(M,g,k)$ which satisfies the dominant energy condition along $\Omega$ (so that $\frak R _{0}\geq 0$). The boundary $\partial\Omega$ is assumed to verify $H^{{\rm APS}}_{0}\geq -2 \widetilde{\lambda}$ (the definition of $\widetilde{\lambda}$ will be given below). Then under the {\rm (APS)} boundary condition, the spectrum of the Dirac--Witten operator on $\Omega$ is an unbounded discrete set of real numbers, such that  any eigenvalue satisfies
\begin{equation}\label{spectreAPS}
	\lambda^{2} \geq \frac{n}{(n-1)}\frak R_{0}  \ .
\end{equation}
\end{thm}

\preuve
The fact that, under the (APS) boundary condition, the spectrum of the Dirac--Witten operator on $\Omega$ is an unbounded discrete set of real numbers, has been proved in the previous lemma.\\
Thanks to the isomorphisms of the several Clifford actions (\eq{clifford}), we can improve (\eq{integrand})
\begin{eqnarray*}
	\omega_{\varphi}(\nu) &=& \left\langle - \widetilde{\D}\varphi+ \frac{1}{2}\big(  \tr_{\ell}\theta + (\tr_{\ell}k)\nu\cdot e_{0}- k^{\partial\Omega}(\nu)\cdot e_{0}\cdot\big)\varphi,\varphi  \right\rangle  \\
	&=& \left\langle - \widetilde{\D}\varphi+ \frac{1}{2}e_{0}\cdot \k_{APS}\cdot \varphi,\varphi  \right\rangle \ ,
\end{eqnarray*}
where we have defined the vector field $\k_{APS}\in \Gamma(TN_{|\partial\Omega})$ as
$$ \k_{APS} = (\tr_{\ell}\theta)e_{0} - (\tr_{\ell}k)\nu + k^{\partial\Omega}(\nu) \ .$$
Let $(\widetilde{\lambda} _{m},\psi_{m})_{m\in \Bbb Z}$  the eigenvalues and eigenspinors of $\widetilde{\D}$ (with the convention that $\widetilde{\lambda}_{m}<0$ if and only if $m<0$). If we suppose that $\Pi_{+}\varphi=0$, then $\varphi=\varphi_{0}+ \underset{m<0}{\sum}\varphi_{m}$, where $\varphi_{m}$ is the $L^{2}$--orthogonal projection of $\varphi$ on the line $\R \psi_{m}$ and $\varphi_{0}\in\Ker \widetilde{\D}$. We can deduce
\begin{equation*}
	\int_{\partial\Omega}\omega_{\varphi}(\nu)= \sum_{m<0} \int_{\partial\Omega}\left\langle (-\widetilde{\lambda}_{m}+ \frac{1}{2}e_{0}\cdot \k_{APS}\cdot) \varphi_{m}, \varphi_{m}\right\rangle \ .
\end{equation*}
It is well known that $\left\langle (-\widetilde{\lambda}_{m}+ \frac{1}{2}e_{0}\cdot \k_{APS}\cdot) \varphi_{m}, \varphi_{m}\right\rangle \geq 0$ if the vector field $(\k_{APS} -2 \widetilde{\lambda}_{m}e_{0})$ is causal and future oriented, for each integer $m <0$. This equivalently means
$$ \tr_{\ell}\theta -2 \widetilde{\lambda}_{m} \geq \left((\tr_{\ell} k)^{2} + \left|k^{\partial\Omega}(\nu)\right|^{2}\right)^{\frac{1}{2}} \ ,$$
for every $m<0$. This condition will be satisfied if 
\begin{equation}\label{bordAPS}
	H^{{\rm APS}}_{0} \geq - 2 \widetilde{\lambda} : =  - 2 \inf_{m<0}\left\{-\widetilde{\lambda}_{m}\right\} \ .
\end{equation}
As a consequence, if one assumes that (\eq{bordAPS}) holds and that $\varphi$ is solution of (APS) for a non--zero $\lambda\in \R$, then (\eq{Int2}) implies
\begin{equation*}
	0 \leq   \left(\frac{n-1}{n}\right) \lambda^{2}\int_{\Omega} \left| \varphi\right|^{2}  -\int_{\Omega} \left\langle \frak R \varphi,\varphi,\right\rangle \ ,
\end{equation*}
which is possible only if 
$$ \left|\lambda\right|^{2} \geq   \frac{n}{(n-1)}\frak R_{0} \ .$$
\qed

\begin{Rq}
	The boundary assumption $H^{{\rm APS}}_{0}\geq -2 \widetilde{\lambda}$ is quite weaker than the one used in \cite{HijMR}, even for the case $k\equiv0$.
\end{Rq}

\begin{Rq}
According to the lower bound of C.~Bär and O.~Hijazi 
\begin{equation}\label{BarHij}
	\widetilde{\lambda} \geq \frac{1}{2}\Vol(\partial\Omega,\ell)^{\frac{-1}{n-1}}\sqrt{\frac{(n-1)}{(n-2)}\mycal{Y}(\partial\Omega,\ell)}  \ ,
\end{equation}
where $\mycal{Y}(\partial\Omega,\ell)$ denotes the Yamabe invariant of the closed manifold $(\partial\Omega,\ell)$, we can replace the boundary assumption $H^{{\rm APS}}_{0}\geq -2 \widetilde{\lambda}$ by the stronger one 
$$ H^{{\rm APS}}_{0}\geq - \Vol(\partial\Omega,\ell)^{\frac{-1}{n-1}}\sqrt{\frac{(n-1)}{(n-2)}\mycal{Y}(\partial\Omega,\ell)} \ .$$
We will use this condition in order to investigate the equality case of (\eq{spectreAPS}).
\end{Rq}

The equality case in (\eq{spectreAPS}) is not easy to treat in general, but we can however deduce some information under the stronger boundary assumption (\eq{BarHij}) introduced in the previous remark.

\begin{thm}\label{APSegal}
Let $\Omega$ be a compact domain of a spacelike spin hypersurface $(M,g,k)$ which satisfies the dominant energy condition along $\Omega$ (so that $\frak R _{0}\geq 0$). The boundary $\partial\Omega$ is assumed to have a positive Yamabe invariant and to verify the boundary inequality $H^{{\rm APS}}_{0}\geq - \Vol(\partial\Omega,\ell)^{\frac{-1}{n-1}}\sqrt{\frac{(n-1)}{(n-2)}\mycal{Y}(\partial\Omega,\ell)}$ . Then equality in (\eq{spectreAPS}) always leads to the existence of an imaginary $\nabla$--Killing spinor on $\Omega$ of Killing number $-\lambda/n$, and the boundary metric $\ell$ is Einstein with positive scalar curvature on $\partial\Omega$. 
\end{thm}

\preuve
Suppose  that equality holds in~(\eq{spectreAPS}). Thereby, there exists a non--zero $\lambda$--eigenspinor field $\varphi$ such that 
\begin{equation}\label{egalAPS}
	P\varphi \equiv 0 \ , \quad \left\langle e_{0}\cdot(\k_{{\rm APS}}-2\widetilde{\lambda}e_{0})\cdot\varphi,\varphi \right\rangle \equiv 0 \ .
\end{equation}
The second equation of~(\eq{egalAPS}) simply says that the vector field $(\k_{{\rm APS}}-2\widetilde{\lambda}e_{0})$ is lightlike, which means in other words
$$H^{{\rm APS}}_{0} = -2\widetilde{\lambda}\ ,$$
but it says that we are in the equality case of (\eq{BarHij}), and so the boundary metric $\ell$ is Einstein with positive scalar curvature on $\partial\Omega$.\\
The first equation of~(\eq{egalAPS}) can be reformulated as 
$$ \nabla_{X}\varphi + \left(\frac{\lambda}{n}\right) X\cdot \varphi =0 \ , \quad\ \forall X \in\Gamma(T\Omega),$$
that is, $\varphi$ is a Killing spinor of Killing number $-\lambda/n$.\qed

\begin{Rq}
	If we replace the boundary assumption of Theorem~\eq{APSegal} by $H^{{\rm APS}}_{0}> -2 \widetilde{\lambda}$ (or by $H^{{\rm APS}}_{0} \geq 0$ as in \cite{HijMR}), then equality in (\eq{spectreAPS}) cannot occur.
\end{Rq}


\subsection{Modified APS Boundary Condition}
The modified Atiyah--Patodi--Singer (mAPS) boundary condition was first introduced in \cite{HijMR}. In this section, our goal is to find  a lower bound for  $\left|\lambda\right|^{2}$ where $\lambda$ is any non--zero complex (a priori) number involving in the following elliptic first order boundary problem
$$(\textrm{mAPS}) \quad \left\{
\begin{array}{ll}
\D \varphi= \lambda\varphi & \textrm{ on } \Omega \\
\Pi_{+}(\varphi+\nu\cdot\varphi)=  0 & \textrm{ on } \partial\Omega
\end{array} \right. \ ,
$$
where $\D$ still denotes the Dirac--Witten operator. As for the (APS) problem, we first prove that the spectrum of (mAPS) is real.
\begin{lem}
The spectrum of {\rm (mAPS)} is a discrete set of real numbers. 
\end{lem}

\preuve Let $\lambda\in \C$ be any eigenvalue of ${\rm (mAPS)}$ with $\varphi$ a corresponding eigenspinor field. We know that 

\begin{eqnarray*}
	(\lambda-\overline{\lambda})\left\|\varphi\right\|_{L^{2}(\Omega)}&=&-\left\langle \nu\cdot\varphi,\varphi\right\rangle_{L^{2}(\partial\Omega)} \\
	&=& -\frac{1}{2}\int_{\partial\Omega}\left\langle \varphi+ \nu\cdot \varphi, \varphi- \nu\cdot\varphi  \right\rangle \\
	&=& 0 \ ,
\end{eqnarray*}
since $\varphi+ \nu\cdot \varphi \in \Sigma^{-}_{|\partial\Omega} $ by definition and  $\varphi- \nu\cdot\varphi  \in \Sigma^{+}_{|\partial\Omega} $. Indeed, we can show that $\Pi_{+}\nu = \nu \Pi_{-}$ (that is a consequence of $\nu\cdot\widetilde{\D}=- \widetilde{\D}\nu\cdot $), and so 
\begin{eqnarray*}
	\Pi_{+}(\nu\cdot\varphi-\varphi) &=& \Pi_{+}(\nu\cdot\varphi+ \nu\cdot\nu\cdot\varphi) \\
	&=& \nu\cdot \Pi_{-} (\varphi+ \nu\cdot\varphi ) \\
	&=& \nu\cdot  (\varphi+ \nu\cdot\varphi )\\
	&=& \nu\cdot\varphi - \varphi \ .
\end{eqnarray*}
\qed

The third main result of this note is a generalisation of the lower bound of \cite{HijMR} for the Dirac--Witten operator under the (mAPS) boundary condition (we recover the inequality of \cite{HijMR} when we set $k\equiv 0$ on $\Omega$).

\begin{thm}\label{mainAPS}
Let $\Omega$ be a compact domain of a spacelike spin hypersurface $(M,g,k)$ which satisfies the dominant energy condition along $\Omega$ (so that $\frak R _{0}\geq 0$). The boundary $\partial\Omega$ is assumed to verify $H^{{\rm APS}}_{0}\geq 0$ . Then under the {\rm (mAPS)} boundary condition, the spectrum of the Dirac--Witten operator on $\Omega$ is an unbounded discrete set of real numbers, such that  any eigenvalue satisfies
\begin{equation}\label{spectremAPS}
	\lambda^{2} \geq \frac{n}{(n-1)}\frak R_{0}  \ .
\end{equation}
\end{thm}

\preuve
The fact that, under the (mAPS) boundary condition, the spectrum of the Dirac--Witten operator on $\Omega$ is an unbounded discrete set of real numbers, has been proved in the previous lemma.\\
We still have 
$$ 	\omega_{\varphi}(\nu) = \left\langle - \widetilde{\D}\varphi+ \frac{1}{2}e_{0}\cdot \k_{APS}\cdot \varphi,\varphi  \right\rangle \ ,$$
where $\k_{APS}$ has been defined in the previous section. If we suppose that ${\Pi_{+}(\varphi+\nu\cdot\varphi)=  0 }$, then  we know (see the proof of the lemma above) that ${\varphi+ \nu\cdot \varphi \in \Sigma^{-}_{|\partial\Omega} }$ and  ${\varphi- \nu\cdot\varphi \in \Sigma^{+}_{|\partial\Omega}}$. Therefore, 
\begin{eqnarray*}
	0 &=& \frac{1}{2}\left\langle\widetilde{\D}(\varphi + \nu\cdot\varphi),\varphi-\nu\cdot\varphi \right\rangle_{L^{2}(\partial\Omega)}\\	
	&=& \left\langle\widetilde{\D}\varphi,\varphi \right\rangle_{L^{2}(\partial\Omega)} \ ,
\end{eqnarray*}
where the second line is obtained thanks to the property $\nu\cdot\widetilde{\D}=- \widetilde{\D}\nu\cdot $. Assume $\varphi$ is solution of (mAPS) for a non--zero $\lambda\in \R$, then (\eq{Int2}) implies
\begin{equation*}
	0 \leq   \left(\frac{n-1}{n}\right) \lambda^{2}\int_{\Omega} \left| \varphi\right|^{2}  -\int_{\Omega} \left\langle \frak R \varphi,\varphi,\right\rangle \ ,
\end{equation*}
since the condition $H^{{\rm APS}}_{0}\geq 0$ implies that $\k_{APS}$ is causal and future oriented. This leads to the desired inequality, namely
$$\lambda^{2} \geq \frac{n}{(n-1)}\frak R_{0} \ .$$\qed

The equality case in (\eq{spectremAPS}) is not easy to treat in general, but we can however deduce some information under a natural additional assumption.

\begin{thm}\label{mAPSegal}
Under the assumptions of Theorem~\eq{mainMIT}, equality in (\eq{spectreMIT}) always leads to the existence of an imaginary $\nabla$--Killing spinor on $\Omega$ of Killing number $-\lambda/n$. If we assume furthermore that ${k^{\partial\Omega}(\nu)=0}$ along the boundary $\partial\Omega$, then  $\partial\Omega$ is an apparent horizon.
\end{thm}

\preuve
Suppose  that equality holds in~(\eq{spectremAPS}). Thereby, there exists a non--zero $\lambda$--eigenspinor field $\varphi$ such that 
\begin{equation}\label{egalmAPS}
	P\varphi \equiv 0 \ , \quad \left\langle e_{0}\cdot\k_{{\rm APS}}\cdot\varphi,\varphi \right\rangle \equiv 0 \ .
\end{equation}
The second equation of~(\eq{egalmAPS}) simply says that the vector field $\k_{{\rm APS}}$ is lightlike, which means in other words
$$H^{{\rm APS}}_{0} = 0 \ ,$$
and the first one reads as 
$$ \nabla_{X}\varphi + \left(\frac{\lambda}{n}\right) X\cdot \varphi =0 \ , \quad\ \forall X \in\Gamma(T\Omega),$$
that is, $\varphi$ is a Killing spinor of Killing number $-\lambda/n$.\\
If we suppose futhermore  ${k^{\partial\Omega}(\nu)=0}$ along the boundary $\partial\Omega$, then $H^{{\rm APS}}_{0} = 0$ can be reformulated as $\tr_{\ell}\theta = \left|\tr_{\ell} k \right|$ which, in the General Relativity literature, is the condition to be an apparent horizon.\qed


\subsection{Chiral Boundary Condition}
In the present section, our aim is to find an estimate for the spectrum of the Dirac--Witten operator under a natural boundary condition associated to a chirality operator. It consists on finding a lower bound for  $\left|\lambda\right|^{2}$ where $\lambda$ is any non--zero complex (a priori) number involving in the following elliptic first order boundary problem
$$(\textrm{CHI}) \quad \left\{
\begin{array}{ll}
\D \varphi= \lambda\varphi & \textrm{ on } \Omega \\
G(\varphi)=  \varphi & \textrm{ on } \partial\Omega
\end{array} \right. \ ,
$$
where $\varphi$ is non--zero (eigen--)spinor field, and where $G\in\End\left(\Sigma_{|\partial\Omega}\right)$ is defined  by the relation  ${G(\psi)= \eps \nu\cdot e_{0}\cdot\psi}$, with $\eps=\pm 1$. The boundary condition $G(\psi)=  \psi$ was originally introduced (as the (APS) condition) to prove some black hole version of the positive mass theorem for asymptotically flat manifolds \cite{GHHP,H1}. Some of our arguments appeared first in \cite{HijZ}, but our conclusion concerning the equality case is much more stronger. As usual in this article, we first prove that the spectrum of (CHI) is real.
\begin{lem}
The spectrum of {\rm (CHI)} is a discrete set of real numbers. 
\end{lem}

\preuve Let $\lambda\in \C$ be any eigenvalue of ${\rm (CHI)}$ with $\varphi$ a corresponding eigenspinor field. We know that 
\begin{eqnarray*}
	(\lambda-\overline{\lambda})\left\|\varphi\right\|_{L^{2}(\Omega)}&=&-\left\langle \nu\cdot\varphi,\varphi\right\rangle_{L^{2}(\partial\Omega)} \\
	&=& -\left\langle \nu\cdot\varphi,e_{0}\cdot e_{0}\cdot\varphi\right\rangle_{L^{2}(\partial\Omega)} \\
	&=&   \left\langle \nu\cdot e_{0} \cdot\varphi, e_{0}\cdot\varphi\right\rangle_{L^{2}(\partial\Omega)} \\
	&=& \eps  \left\langle \varphi, e_{0}\cdot\varphi\right\rangle_{L^{2}(\partial\Omega)} \ , 
\end{eqnarray*}
which is real since $e_{0}\cdot$ is Hermitian with respect to the scalar product $\left\langle *, * \right\rangle$. \qed

For later use we introduce a geometric quantity.
\begin{déf} We set
$$H_{0}^{{\rm CHI}}:= \inf_{\partial\Omega}\left\{\tr_{\ell}\theta + \eps (\tr_{\ell}k)  \right\}  \ .$$
\end{déf}

The last main result of this note is a generalisation of the lower bound of \cite{HijMR} for the Dirac--Witten operator under the (CHI) boundary condition (we recover the inequality of \cite{HijMR} when we set $k\equiv 0$ on $\Omega$).

\begin{thm}\label{mainCHI}
Let $\Omega$ be a compact domain of a spacelike spin hypersurface $(M,g,k)$ which satisfies the dominant energy condition along $\Omega$ (so that $\frak R _{0}\geq0$). The boundary $\partial\Omega$ is assumed to verify $H^{{\rm CHI}}_{0} \geq 0$. Then under the {\rm (CHI)} boundary condition, the spectrum of the Dirac--Witten operator on $\Omega$ is an unbounded discrete set of real numbers, such that  any eigenvalue satisfies
\begin{equation}\label{spectreCHI}
	\left|\lambda\right|^{2} \geq \frac{n}{(n-1)}\frak R_{0}  \ .
\end{equation}
\end{thm}

\preuve
The fact that, under the (CHI) boundary condition, the spectrum of the Dirac--Witten operator on $\Omega$ is an unbounded discrete set of real numbers, has been proved in the previous lemma.\\
We  now need to give some elementary properties of the boundary spinorial endomorphism $G$ (the proof is left to the reader).
\begin{prop}\label{symCHI}
The endomorphism $G$ is symmetric, isometric with respect to $\left\langle\ast,\ast\right\rangle$, anticommutes to the action of $\nu\cdot$  and $e_{0}\cdot$, and commutes to each $e_{k}\cdot,\, (k \geq 2)$.
\end{prop}
The important fact is that the boundary condition of (CHI) allows us to express the boundary integrand of (\eq{Int5}) in terms of the extrinsic curvature tensors $k$ and $\theta$. 

\begin{lem}\label{G}
If $G(\varphi)=\varphi$ then along the boundary $\partial\Omega$ we have
\begin{equation}\label{integrandCHI}
	\omega_{\varphi}(\nu)=\frac{1}{2}\big(  \tr_{\ell}\theta  +\eps (\tr_{\ell} k) \big)  \left|\varphi\right|^{2}   \ .
\end{equation}
\end{lem}
\noindent
\textsc{Proof of Lemma~\eq{G}:}
We use Equation~(\eq{integrand}) 
\begin{eqnarray*}
	\omega_{\varphi}(\nu) &=& \left\langle \nu\cdot  e_{j}\cdot \widetilde{\nabla}_{e_{j}}\varphi + \frac{1}{2}\Big( (\tr_{g} k ) \nu\cdot e_{0}\cdot - k(\nu)\cdot e_{0}\cdot + (\tr_{\ell}\theta)  \Big) \varphi,\varphi\right\rangle\\
	&=& \left\langle \nu\cdot  e_{j}\cdot \widetilde{\nabla}_{e_{j}}\varphi + \frac{1}{2}\Big( (\tr_{\ell} k ) \nu\cdot e_{0}\cdot - k^{\partial\Omega}(\nu)\cdot e_{0}\cdot + (\tr_{\ell}\theta)  \Big) \varphi,\varphi\right\rangle \\
	&=& \frac{1}{2}\big(  \tr_{\ell}\theta  +\eps (\tr_{\ell} k) \big) \left|\varphi\right|^{2}  \ .
\end{eqnarray*}
when we take Proposition~\eq{symCHI}  and  the assumption $G(\varphi)=\varphi$ into account. Assume $\varphi$ is solution of (CHI) for a non--zero $\lambda\in \R$, then (\eq{Int2}) implies
\begin{equation*}
	0 \leq   \left(\frac{n-1}{n}\right) \lambda^{2}\int_{\Omega} \left| \varphi\right|^{2}  -\int_{\Omega} \left\langle \frak R \varphi,\varphi,\right\rangle \ ,
\end{equation*}
since the condition $H^{{\rm CHI}}_{0}\geq 0$ induces the non--negativity of the boundary integral in  (\eq{Int2}).  This leads to the desired inequality, namely
$$\lambda^{2} \geq \frac{n}{(n-1)}\frak R_{0} \ .$$\qed

On the contrary to the other boundary conditions, equality case in (\eq{spectreCHI})  can be  treated in general. 

\begin{thm}\label{CHIegal}
Under the assumptions of Theorem~\eq{mainCHI}, equality in (\eq{spectreCHI}) always leads to the existence of an imaginary $\nabla$--Killing spinor on $\Omega$ of Killing number $-\lambda/n$, and the second fundamental form of  $\partial\Omega$ is proportional to $k^{\partial\Omega}$, since
$$\theta +\eps  k^{\partial\Omega} \equiv 0 \ .$$
In particular, $\partial\Omega$ is an apparent horizon.
\end{thm}

\preuve
Suppose now that equality holds in~(\eq{spectreCHI}). Thereby, there exists a non--zero $\lambda$--eigenspinor field $\varphi$ such that 
\begin{equation}\label{egalCHI}
	P\varphi \equiv 0 \ , \quad H^{{\rm CHI}}_{0} \equiv 0 \ .
\end{equation}
The second equation of~(\eq{egalCHI}) simply says that  $\tr_{\ell}\theta + \eps (\tr_{\ell} k)= 0 $ on $\partial \Omega $. The first equation of~(\eq{egalCHI}) can be reformulated as 
$$ \nabla_{X}\varphi + \left(\frac{\lambda}{n}\right) X\cdot \varphi =0 \ , \quad\ \forall X \in\Gamma(T\Omega),$$
that is, $\varphi$ is a Killing spinor.  We consider now any vector field $X\in\Gamma(T\partial\Omega)$ tangent to the boundary, and we compute
\begin{eqnarray*}
	 -\left(\frac{\lambda}{n}\right) X \cdot \varphi &=& \nabla_{X}\varphi\\
	&=&  \eps \nabla_{X}(\nu\cdot e_{0} \cdot \varphi) \\
	&=&  \eps\big((\nabla_{X}\nu)\cdot e_{0} \cdot\varphi + \nu\cdot (\nabla_{X}e_{0}) \cdot \varphi + \nu\cdot e_{0} \cdot \nabla_{X}\varphi\big)\\
	&=&  \eps \left(  \overline{\nabla}_{X}\nu\cdot e_{0}\cdot -k(X,\nu) - \nu\cdot k(X)\cdot - \left(\frac{\lambda}{n}\right)\nu\cdot e_{0}\cdot X\cdot  \right)\varphi \\
	&=&  \eps\big( -\theta(X)\cdot e_{0} \cdot  -k(X,\nu) - \nu\cdot k(X)\cdot\big)\varphi -  \left(\frac{\lambda}{n}\right)X\cdot\varphi \ .
\end{eqnarray*}
This finally leads to (after a left multiplication by $e_{0}\cdot$)
\begin{equation}\label{meanCHI1}
	\forall X\in\Gamma(T\partial\Omega) \ , \quad \left( \theta(X)\cdot -k(X,\nu)e_{0}\cdot  + \nu\cdot e_{0}\cdot k(X) \cdot  \right)\varphi=0 \ .
\end{equation}
The condition~(\eq{meanCHI1}) can be improved, using once again $G(\varphi)=\varphi$
\begin{eqnarray*}
	0&=& \left( \theta(X)\cdot -k(X,\nu)e_{0}\cdot  + \nu\cdot e_{0}\cdot k(X) \cdot  \right)\varphi \\
	&=& \left( \theta(X)\cdot -k(X,\nu)e_{0}\cdot  + \nu\cdot e_{0}\cdot k^{\partial\Omega}(X) \cdot  + \nu\cdot e_{0}\cdot k(X,\nu) \nu\cdot \right)\varphi \\
	&=& \left( \theta(X)\cdot   + \nu\cdot e_{0}\cdot k^{\partial\Omega}(X) \cdot   \right)\varphi \\
	&=& \left( \theta(X) +\eps  k^{\partial\Omega}(X)   \right) \cdot \varphi \ ,
\end{eqnarray*}
for every $X$  tangent to the boundary. This induces  $\theta +\eps  k^{\partial\Omega} \equiv 0$  on $\partial\Omega$. 
 \qed

\begin{Rq}
	The Gauss equation explicitely gives the Ricci curvature of $g$ restricted to the boundary when equality is achieved in~(\eq{spectreCHI})
\begin{eqnarray*}
	\Ric^{g,\partial\Omega} &=& \frac{\Scal^{\ell}}{(n-1)}\ell - (\tr_{\ell}\theta)\theta + \theta\circ\theta  \\
	&=& \frac{\mycal{Y}(\partial\Omega,\ell)}{(n-1)}\ell - (\tr_{\ell}k)k^{\partial\Omega} + k^{\partial\Omega}\circ k^{\partial\Omega} \ , 
\end{eqnarray*}
since $\ell$ is an Eintein Yamabe metric, and  $\theta +\eps  k^{\partial\Omega} \equiv 0$  on $\partial\Omega$ $(\eps^{2}=1)$. 
\end{Rq}


\end{document}